\input amstex

\documentstyle{amsppt}
\document
\topmatter

\define\AC{\Cal AC^{\perp}}

\title
On some $\AC$ manifolds
\endtitle

\author
W\l odzimierz Jelonek
\endauthor

\par
\medskip
\abstract{We give a description of all Gray $\AC$ manifolds $(M,g)$  whose
Ricci tensor has two eigenvalues of multiplicity $1$ and $\dim M-1$.}
\endabstract
\thanks{key words and phrases : Gray manifold, Einstein-Weyl manifold, MSC: 53C25}\endthanks
\endtopmatter

\define\G{\Gamma}

\define\e{\epsilon}
\define\n{\nabla}
\define\om{\omega}
\define\w{\wedge}
\define\k{\diamondsuit}

\define\p{\parallel}
\define\a{\alpha}

\define\lb{\lambda}
\define\Lb{\Lambda}
\define\A{\Cal A}
\define\1{D_{\lb}}
\define\2{D_{\mu}}

\define\bt{\beta}

\define\De{\Cal D}

\define\de{\delta}
\define\r{\rightarrow}
\par
\medskip
{\bf 0. Introduction.} In this paper we study complete Riemannian
manifolds whose Ricci tensor $\rho$ satisfies the relation
$$\n_X\rho(X,X) =\frac{2X \text{Scal}}{n+2} g(X,X)\tag 0.1 $$
where $ {\text{Scal}}= \text{tr}_g \rho$ is the scalar curvature of $(M,g)$ and $n=\dim M$.
This  property was studied by A. Gray  [Gr] and by A. Besse [B].
A. Gray called Riemannian manifolds satisfying $(0.1)$  $\AC$
manifolds;  we shall also call them  Gray manifolds. We describe all simply connected and real analytic Gray $\AC$ manifolds whose
Ricci tensor has two eigenvalues $\lb,\mu$  of multiplicity $\dim M-1$ and $1$ .

If $([g],D)$ is an Einstein-Weyl structure on a compact manifold $M$, and $g_0$ is the Gauduchon metric  in the conformal class $[g]$,  then   $(M,g_0)$ is an $\Cal{AC}^{\perp}$ manifold (see [J-1]) with the above property. If $\om_0$ is the Lee form of $g_0$, i.e. $Dg_0=\om_0\otimes g_0$, then the Weyl structure determined by the pair $(g_0,-\om_0)$ is also Einstein-Weyl (see [J-1]). Hence if $([g],D)$ is an Einstein-Weyl structure on a compact manifold $M$,  and $g_0$ is the Gauduchon metric of $([g],D)$, then on $M$ there is  another Einstein-Weyl structure $([g],D^-)$ such that $D^-g_0=-\om_0\otimes g_0$.  Thus $(M,g_0)$ is then a Riemannian manifold which admits a pair of Einstein-Weyl structures determined by the pairs $(g_0,\om_0)$ and $(g_0,-\om_0)$.  More generally a Riemannian manifold $(M,g)$ which admits two Einstein-Weyl structures determined by pairs $(g,\om)$ and $(g,-\om)$ is an $\Cal{AC}^{\perp}$ manifold whose
Ricci tensor has two eigenvalues $\lb,\mu$  of multiplicity  $\dim M-1$ and $1$.
We prove that Gray manifolds with two eigenvalues $\lb,\mu$  of multiplicity $\dim M-1$ and $1$ always admit a conformal vector field $\xi$ which is an eigenfield of the Ricci tensor of $(M,g)$. The field $\xi$ is a Killing vector field or a conformal closed gradient vector field.

In particular we give a  description of all complete, real analytic, simply connected Riemannian manifolds $(M,g)$ admitting a pair of Einstein Weyl structures $(g,\om),(g,-\om)$. We show that on such manifolds $\om(X)=g(\xi,X)$ where $\xi$ is either Killing vector field or  a closed gradient conformal field. In the second case we give a complete classification of such manifolds. We prove that in this case  the two Einstein-Weyl structures are conformally Einstein.

For  general $\Cal{AC}^{\perp}$ manifolds we prove that the eigenvalues $\lb,\mu$ of the Ricci tensor satisfy the equations: $ n\mu-2(n-1)\lb=C_0=const$ in the case of a non-Killing conformal vector field and $(n-4)\lb+2\mu=C_0=const$ in the case of a Killing vector field.

 We also give new examples of compact $\Cal{AC}^{\perp}$ manifolds diffeomorphic to the sphere $S^{n},n>2$.

\par
\bigskip
{\bf 1. $\AC$ manifolds and Einstein-Weyl manifolds.}

 By an $\AC$ {\it  manifold }(see [Gr],[B]) we mean a Riemannian manifold $(M,g)$ satisfying the condition
$$\frak C_{X Y Z}\n_X\rho(Y,Z)=\frac 2{\dim M+2}\frak C_{X Y Z}X \text{Scal} g(Y,Z),  \tag 1.1$$
where $\rho$ is the Ricci tensor of $(M,g)$,  $\text{Scal}$ is the scalar curvature  and $\frak C$ means
the cyclic sum.

In [J-1] it is proved that a Riemannian manifold
$(M,g)$ is an $\AC$ manifold if and only if the Ricci endomorphism
$Ric$ of $(M,g)$ is of the form $Ric=S+\frac{2}{n+2}  \text{Scal} Id$ where
$S$ is a Killing tensor, and
$n=$dim$M$. Let us recall that a $(1,1)$ tensor $S$ on a Riemannian
manifold $(M,g)$ is called a {\it Killing tensor} if $g(\n S(X,X),X)=0$
for all $X\in TM$.

 Define  the integer-valued function  $E_S(x)=($ the number of distinct eigenvalues of $S_x$) and
 set $M_S=\{ x\in M:E_S$ is constant in a neighborhood of $x\}$.
 The set $M_S$ is open and dense in $M$ and the eigenvalues $\lb_i$ of $S$ are distinct and
 smooth in each component $U$ of $M_S$.  Let us denote by $D_{\lb_i}$
 the eigendistribution corresponding to $\lb_i$. We have  (see [J-1])
\medskip
{\bf Proposition.} {\it Let } $S$ {\it be a Killing  tensor on} $M$ {\it and U be a component of }
 $M_S$ {\it and} $\lb_1,\ldots,\lb_k \in C^{\infty}(U)$
{\it be eigenfunctions of } $S$. {\it Then for all} $X\in D_{\lb_i}$ {\it we have}
   $$ \n S(X,X)=-\frac12 \n \lb_i\| X \| ^2\tag 1.2$$
{\it and}  $D_{\lb_i}\subset \ker d\lb_i$. { \it If } $i\ne j$ {\it and} $X\in \G(D_{\lb_i}),Y \in \G(D_{\lb_j})$  {\it then}
$$ g(\n_X X ,Y)=\frac12\frac{Y\lb_i}{\lb_j-\lb_i}\| X \| ^2.\tag 1.3$$

If $T(X,Y)=g(SX,Y)$ is a Killing tensor on $(M,g)$ and $c$ is a geodesic on $M$ then the function  $\phi(t)=T(\dot c(t),\dot c(t))$ is constant on the domain of $c$. In fact $\phi'(t)=\n_{\dot c(t)}T((\dot c(t),\dot c(t))=0$.

A {\it conformal vector field} $\xi$ on a Riemannian manifold satisfies
the relation $L_{\xi}g=\a g$.   This is equivalent to
$\n_X\om(Y)+\n_Y\om(X)=\a g(X,Y)$ where $\om=g(\xi,.)$.

Finally let us recall that a Riemannian manifold with Killing Ricci tensor is called after A. Gray
an $\A$ {\it manifold} ( it is an $\AC$ manifold with constant scalar curvature).

 We start with some basic facts
  concerning Einstein-Weyl geometry. For more details see [T], [PS1],
[PS2].

Let $M$ be an $n$-dimensional manifold with a conformal structure
$[g]$ and a torsion-free affine
 connection $D$. This defines an {\it Einstein-Weyl} (E-W) {\it structure}
if  $D$ preserves the conformal structure, i.e. there exists a
1-form $\om$ on $M$ such that
$$Dg=\om\otimes g \tag 1.4$$
and the Ricci tensor $\rho^D(X,Y)=\text{tr}\{Z\r R^D(Z,X)Y\}$ of $D$, where $R^D$ is the curvature tensor of $(M,D)$, satisfies the condition
$$\rho^D(X,Y)+\rho^D(Y,X)=\bar\Lb g(X,Y) \text{  for every  } X,Y\in TM$$
for some function $\bar\Lb\in C^{\infty}(M)$. P.  Gauduchon ([G]) proved
 the fundamental theorem that if $M$ is compact then there
exists a Riemannian metric $g_0\in [g]$ for which the Lee form $\om_0$ associated to $g_0$ is co-closed-$\delta \om_0=0$
where $Dg_0=\om_0\otimes g_0$, and $g_0$ is unique up to homothety. The metric $g_0$ is  called {\it Gauduchon} or the {\it standard} metric of the E-W structure $([g],D)$ on $M$.
 Note that if $M$ is not compact, a standard metric may not exist.

  Let $\rho$ be the Ricci tensor of $(M,g)$
  and let us denote by $Ric$ the Ricci
 endomorphism of $(M,g)$, i.e. $\rho(X,Y)=g(X, Ric Y)$.
We recall an important theorem (see [PS1]):
\par
\medskip
{\bf Theorem.  }  {\it A metric g and a 1-form $\om$ determine an
E-W structure if and only if there exists a function $\Lb\in
C^{\infty}(M)$ such that
$$\rho+\frac14(n-2) \Cal D\om=\Lb g  \tag 1.5$$
where $\Cal D\om(X,Y)=(\n_X\om)Y+(\n_Y\om)X+\om(X)\om(Y)$ and
$n=\dim \ M$. If (1.5) holds then}
$$\bar\Lb=2\Lb+\text{div}\om-\frac12 (n-2)\p\om^{\sharp}\p^2\tag 1.6$$

\bigskip
 Compact Einstein-Weyl manifolds with the
Gauduchon metric are Gray manifolds. By a {\it  Riemannian manifold with a pair of Einstein-Weyl structures} we mean a manifold $(M,g)$ for which there exists a 1-form $\om$ such that $(g,\om),(g,-\om)$ determine two Einstein Weyl structures  $([g],D),([g],D^-)$ with $Dg=\om\otimes g, D^-g=-\om\otimes g$. Note that every  Einstein-Weyl structure $(D,[g])$ on a compact manifold $M$ determines a Riemannian manifold $(M,g)$ with two Einstein-Weyl structures, namely $([g], D),([g], D^-)$ where $g$ is the Gauduchon metric for the conformal manifold $(M,[g])$  and $D^-g=-\om_0\otimes g$ where $\om_0$ is the Lee form of $g$ with respect to $D$ (see [J-1]).

It is not difficult to prove that any Riemannian manifold admitting a pair of E-W structures is an $\Cal {AC}^{\perp}$ manifold with two eigenvalues $\lb,\mu$ of the Ricci tensor of multiplicities $n-1,1$ and such that $\lb\ge\mu$ (see [Go]). In this case
$\om(X)=g(\xi,X)$ where $\xi$ is a conformal vector field.  We will prove  the converse: if $(M,g)$ is a simply connected, real analytic $\Cal{AC}^{\perp}$ manifold with two eigenvalues $\lb,\mu$ of multiplicities $n-1,1$ and $\lb>\mu$ at least at one point then $(M,g)$ admits a pair of  Einstein-Weyl structures $( g, \om)$ and $( g,-\om)$.

\bigskip
{\bf 2. $\AC$ manifolds $(M,g)$ whose Ricci tensor has two eigenvalues of
multiplicity $\dim M-1$ and $1$.}  We start with:
\bigskip
{\bf Theorem 1. }{\it  Let $(M,g)$ be a real analytic complete,
simply connected $\AC$ manifold whose Ricci tensor has two
eigenvalues $\lb,\mu$ of multiplicity $n-1,1$ in the set
$V=\{x\in M:\lb(x)\ne\mu(x)\}$.  Then the set $N=\{x\in
M:\lb(x)=\mu(x)\}$ has an empty interior and
there exists a conformal vector field $\xi$ such that
$g(\xi,\xi)=|\lb-\mu|$ and $(Ric-\mu Id)\xi=0$. If $\xi$ is a
Killing vector field then $\lb'$ is constant while if $\xi$ is conformal not
Killing then $\mu'$ is constant. Here  $\lb',\mu'$ are eigenvalues
of the Killing tensor $S=Ric-\frac{2\text{Scal}}{n+2 }Id$ corresponding to
$\lb,\mu$ respectively.}

\medskip

{\it Proof.}  Let  $V_+=\{x\in
M:\lb(x)>\mu(x)\}$,  $V_-=\{x\in M:\lb(x)<\mu(x)\}$. Note that
$M=V_+\cup V_-\cup N$. Let
$\De$ be the distribution on $V$ defined by $\2=\De=\ker(Ric-\mu I)$.
Let us define on $V$ the tensor $m=g(p_{\De}.,p_{\De}.)$ where
$p_{\De}:TM\r \De$ is the orthogonal projection onto $\De$.   Then
$\rho=\lb g+(\mu-\lb)m$. Let us assume that $V_+\ne\emptyset$.  The
distribution $\De_{|V_+}$ is locally spanned by a unit vector
field $\xi_0$.  We have $Ric\xi_0=\mu\xi_0$.  Hence on $V_+$,
$$ \rho=\lb g- \om\otimes\om$$ where $\om=g(\xi, .)$ and
$\xi=\sqrt{\lb-\mu}\xi_0$.

Now we shall check when the tensor
$S=Ric-\frac2{n+2}\text{Scal} g$ is a Killing tensor. Note that for $T(X,Y)=g(SX,Y)$ we have
$T=\lb'g-\om\otimes\om$ where $\lb'=\lb-\frac2{n+2}\text{Scal}$.   Hence
$$\n_XT(X,X)=X\lb'g(X,X)-2\n_X\om(X)\om(X)=X\lb'g(X,X)-h(X,X)\om(X)$$
where $h$ is the symmetric $(2,0)$ tensor defined by
$h(X,Y)=\n_X\om(Y)+\n_Y\om(X)$.  Consequently, $T$ is a Killing tensor if
and only if $$X\lb'g(X,X)-h(X,X)\om(X)=0.\tag 2.1 $$ If
$\om(X)=0$ then from (2.1) it follows that $d\lb'(X)=0$.   Thus
$d\lb'=\a\om$ for some $\a\in C^{\infty}(V_+)$.
Consequently
$$ \om(X)(\a g(X,X)-h(X,X))=0$$  for every $X\in TM_{|V_+}$.

Let $F(X,Y)=\a g(X,Y)-h(X,Y)$.  Then $F$ is a symmetric $(2,0)$
tensor and $$\om(X)F(X,X)=0\tag 2.2$$ for every $X\in TM_{|V_+}$.
Let us write $X=tY+\xi$ where $Y\in\1$, $\om(Y)=0$,  $t\in\Bbb R$. Then
$\om(X)=g(\xi,\xi)\ne0$ and from (2.2) we get
$F(tY+\xi,tY+\xi)=0$.  Consequently,
$t^2F(Y,Y)+2tF(Y,\xi)+F(\xi,\xi)=0$.  Since $t\in\Bbb R$ is
arbitrary, we obtain $F(Y,Y)=F(\xi,Y)=F(\xi,\xi)=0$.  Thus $F=0$.
Hence $h(X,Y)=\a g(X,Y)$ and $\xi$ is a local conformal field.
Since $h(X,Y)=\n_X\om(Y)+\n_Y\om(X)$ we get $n\a=-2\de\om$.  We
also have $$d\lb'=\a\om\tag 2.3$$ and $0=d\a\w\om+\a d\om$. Hence
in the set $U_+=\{x\in V_+:\a(x)\ne0\}$ we obtain
$d\om=-d\ln|\a|\w\om$. It follows that in $U_+$ the distribution
$\1=\ker\om$ is integrable.  Thus for $X,Y\in \1$ we get
$\n_X\om(Y)=g(\n_X\xi,Y)=-g(\xi,\n_XY)=-g(\xi,\n_YX)=g(\n_Y\xi,X)=\n_Y\om(X)$.
Consequently  $\n_X\om(Y)=\frac12\a g(X,Y)$ and
$p_{\1}(\n_X\xi)=\frac12\a X$.     Hence for a vector field $\xi$
and $X\in\1$ we get $\n_X\xi=\frac{\alpha}2X+\phi(X)\xi$ for
$\phi=\frac12d\ln g(\xi,\xi)$. Now  for $X,Y\in\1$ we have

$$\gather 2\n_X\n_Y\xi=X\alpha Y+\alpha\n_XY+2X\phi(Y)\xi+\phi(Y)(\alpha
X+2\phi(X)\xi),\\
 2\n_Y\n_X\xi=Y\alpha X+\alpha\n_YX+2Y\phi(X)\xi+\phi(X)(\alpha
 Y+2\phi(Y)\xi),\\
2\n_{[X,Y]}\xi=\alpha[X,Y]+2\phi([X,Y])\xi\endgather$$
Thus $2R(X,Y)\xi=X\alpha Y-Y\alpha X-\phi(X)\alpha Y+\phi(Y)\alpha
X $ for $X,Y\in\ker\om$. Let $\bt=d\alpha-\alpha\phi$.
Consequently, taking the trace over vectors perpendicular to $\xi$ in the previous formula we get
$0=\rho(Y,\xi)=\frac12(\bt(Y)-(n-1)\bt(Y))=-\frac{n-2}2\bt(Y)$. Hence
$\bt(X)=0$ for $X\in\1$ and
$$R(X,Y)\xi=0$$ for $X,Y\in\1$.
We have $ d\ln|\alpha|(X)=\phi(X)$ for $X\in \1$. Hence for $X\in
\1$, $\n_X\xi=\frac{\alpha}2X+d\ln|\alpha|(X)\xi$ and
$\n_X\om(\xi)=d\ln|\alpha|(X)\om(\xi)$. Thus
$\n_{\xi}\om(X)=-d\ln|\alpha|(X)\om(\xi)$ and
$d\om(X,\xi)=2d\ln|\alpha|(X)\om(\xi)$.  On the other hand,
$d\om=-d\ln|\alpha|\w \om$ and $d\om(X,\xi)=-d\ln|\alpha|(X)\om(\xi)$.
Consequently $d\ln|\alpha|(X)=0$ for $X\in\1$. It follows that $d\om=0$ and
 $\n_X\xi=\frac12\alpha X$ for all $X\in TM$. Thus  $ R(X,Y)\xi=\frac12 d\a(X)Y-\frac12d\a(Y)X$, and consequently
$ \rho(Y,\xi)=-\frac{n-1}2d\a(Y)$.  It follows that $d\a=-\frac{2\mu}{n-1}\om$.

 We also have $d
g(\xi,\xi)(X)=2g(\n_X\xi,\xi)=\alpha\om(X)=d\lb'(X)$. Note that
$d\mu'=0$ in $M-N$ where $\mu'$ is an eigenvalue of the Killing tensor
associated with $Ric$ corresponding to the 1-dimensional
distribution $\2$  (see (1.3)).  Note that $g(\xi,\xi)=\lb'-\mu'$.

Analogously we construct a conformal local vector field $\xi$ in $V_-$. In particular in $V_-$  we have $d\lb'=-\alpha\om$ where $n\alpha=-2\delta \om$ and $\n_X\xi=\frac{\alpha}2X$. Note that if $\xi,\xi'$ are two local fields constructed as above then $\xi=\pm\xi'$ in the intersection of their domains. Now let us define $M_-=\text{int}\{x\in M-N:\a(x)=0\}$ and $M_+=\{x\in M-N:\a(x)\ne0\}$.

We now show that if $M_-\ne\emptyset$ then  $M_+=\emptyset$ hence $\xi$ is a Killing vector field  (see [J-1]).  Let $x_0\in M_-$ and $x_1\in M_+$.  If there exists a smooth curve $c$ joining  $x_0$ and $x_1$ and such that $\text{im }c\cap N=\emptyset$ then by  analytic continuation we see that $\lb'$ is constant along this curve i.e.  $\lb'$ is constant on a certain neighborhood $U_x$ of any $x\in \text{im} c$.  Analogously $\mu'$ is constant and hence $\a(x_1)=0$, a contradiction.  A curve as above indeed exists: if $U$ is a connected component of $x_1$ in $M_+$ then there exists a geodesic $c$   with $c(0)=x_0 $ and  $c(1)\in U$ which doesn't meet $N$  (see  [J-1]).

If  $N=\emptyset$  then $\xi$ exists and is smooth in a neighborhood of any point $x\in M$. Note that in the case of a Killing vector field the constructed vector field is smooth also if $N\ne\emptyset$, and then one of the sets $V_+,V_-$ vanishes. Next we show that a similar situation holds for a conformal vector field.

Let us assume now that $\text {int}\{x:\a(x)=0\}=\emptyset$. In every component of the set $M-N$ and in the interior of $N$ the function $\mu'$ is constant. Since $\mu'$ is continuous it is constant on $M$. The set $\{x:\a(x)\ne0\}$ is dense in $M-N$. In this set $\n_X\xi=\frac\a2X$.  If $x_0\in M- N$ then in a neighborhood of $x_0$ the local field $\xi$ is smooth and $\n_X\xi=\frac\a2X$ also if $\a(x_0)=0$ (we take the limit). Hence in $M-N$ we have $d\a=-\frac{2\mu}{n-1}\om$.  Note that  the function $\a^2$ is independent of the choice of the local field $\xi$, and $\a^2=\frac{||\n\lb'||^2}{|\lb-\mu|}$, hence $\a$ is continuous in $M-N$.  Note that $\a d\a=-\frac{2\mu}{n-1}\a\om=-\frac{2\mu}{n-1}d\lb'$. On the other hand, since $\mu'$ is constant it follows that $n\mu-2(n-1)\lb=C_0=const$ and $d\lb=\frac n{2(n-1)}d\mu$.   Since $\lb-\mu=\lb'-\mu'$ and $\mu'$ is constant, we get $d\lb'=d\lb-d\mu$ and   $(n-1) \a d\a=-2\mu ( d\lb- d\mu)=\frac{2-n}{n-1}\mu d\mu$.   Hence in any component of $M-N$ we have   $\pm (n-1)^2\a^2+(n-2)\mu^2=C_1=const$ (the minus sign is in the components contained in $V_-$).  Since $\mu'$ is constant it follows that $\lb',\lb,\mu$ are real analytic functions on the whole of $M$   (we have  $(n-1)\lb+\mu=\text{Scal}, n\mu-2(n-1)\lb=C_0=const$   and the scalar curvature $\text{Scal}$ of $(M,g)$ is real analytic).

 Now we show that   $\n\lb'=0$ in $N$.  It is clear that $\n\lb'=0$ in the interior of $N$. Let $x_0\in \partial N$.  Then there exists a sequence of points $x_m$ from one component of $M-N$  such that  $lim x_m=x_0$. This component for example is contained in $V_+$.  In that component  $(n-1)^2\a^2+(n-2)\mu^2=C_0=const$.  In particular  $lim_{m\rightarrow\infty}\a^2(x_m)=\frac1{(n-1)^2}(C_0-(n-2)\mu(x_0)^2)$ is finite which means that $\n\lb'=0$ at $x_0$.

Now let assume that $N\ne\emptyset$ and let $x_0\in N$.  Let $c(t)$ be a geodesic with $c(0)=x_0,\|\dot c(0)\|=1$.  If $S=Ric-\frac{2\text{Scal}}{n+2}g$ then $T(.,.)=g(S.,.)$ is a Killing tensor and   $T-\lb'g=\pm\om\otimes\om$.  Note that $\lim_{t\rightarrow 0}T(\dot c(t),\dot c(t))=\mu'$ and since  $T(\dot c(t),\dot c(t))$ is constant it follows that  $g(\xi(c(t)),\dot c(t))^2=\pm(\lb'-\mu')$.  Hence $g(\frac1{\sqrt{|\lb'-\mu'|}}\xi(c(t)),\dot c(t))=\pm 1$.  Since both vectors have length 1 it follows that $\xi(c(t))=\pm\sqrt{|\lb'-\mu'|}\dot c(t)$.  If we introduce  geodesic polar normal coordinates centered at $x_0$  then   we can assume that $\xi=\sqrt{|\lb'-\mu'|}\frac{\partial}{\partial t}$ where $t$ is a radial coordinate.  In particular  $\n\lb'=\a\sqrt{|\lb-\mu|}\frac{\partial}{\partial t}$ on $M-N$ and $\n\lb'=0$ on $N$.  Hence $d\lb'(X)=0$ for any $X$ tangent to the geodesic sphere and $\lb'$ is constant on geodesic spheres with center $x_0$. Thus $\lb'$ is a function of the radial coordinate $t$ only.

If there existed a hypersurface $\Sigma\subset B(x_0,r)$ with $x_0\in\Sigma$ such that the intersection of $\Sigma$ with every geodesic sphere $S(x_0,s),s<r,$ is not empty,  and on which $\lb'-\mu'=0$, then we would have $\lb'=\mu'$ on the whole of $B$.
In fact, if for $x_1\ne x_0$ we have $\lb'(x_1)=\mu'$ then  $\lb'(x)=\mu'$ on the sphere $S(x_0,s)$ containing the point $x_1$.  Hence the sphere $S(x_0,s)$ would be a hypersurface $\Sigma$ in a ball $B(x_1,t)$, where $t$ is the injectivity radius at the point $x_1$, which has the above property. It follows that the set $N$ has a non-empty interior $U$. Since $\lb,\mu $ are real analytic and $\lb=\mu$ on $U$ it follows that $\lb=\mu$ on $M$, a contradiction.   It follows that in a geodesic ball $B(x_0,r)$ where $r$ is the injectivity radius of $M$ at $x_0$, the set $\{x\in B(x_0,r): \lb'(x)=\mu'\}$ consists of only the point $x_0$. Hence the set $N$ consists of isolated points. This implies that $M-N$ has only one component.
It follows that $\lb'-\mu'$ has a constant sign and one of the sets $V_+,V_-$ is empty.

For $x_0\in N$ let us define a local field $\xi$ on the geodesic ball  $B(x_0,r)$ where $r$ is the injectivity radius at $x_0$ by $\xi=\sqrt{|\lb'-\mu'|}\frac{\partial}{\partial t}$ where $t$ is the radial coordinate.  Then since $M$ is simply connected there exists  a continuous vector field $\xi$ on $M$ such that $\rho=\lb g -\om\otimes\om$ if $V_-=\emptyset$ or  $\rho=\lb g +\om\otimes\om$ if $V_+=\emptyset$ and $\om=g(\xi,.)$ (we take the germs of local fields $\xi$ as points of a two-sheeted covering of $M$). The field $\xi$ is smooth on $M-N$ and defines a function $\a$ on $M-N$. The function $\a$ is smooth on $M-N$. We shall show that $\xi$ is in fact smooth.

First we show that the function $\a$ has a continuous extension to the points $x_0\in N$. We assume for example that $V_-=\emptyset$. Then the function $\a^2$ is real analytic and $(n-1)^2\a^2+(n-2)\mu^2=C_1=const$ in $M-N$ hence $\a^2$ has a real analytic extension on the whole of $M$.  If $\a(x_0)^2>0$ then  $\a$ has a constant sign in a certain neighborhood of $x_0$ and $\a=\pm\frac1{n-1}\sqrt{C_1-(n-2)\mu^2}$ is a smooth extension of $\a$ on the neighborhood of $x_0$. If an extension of $\a^2$ satisfies $\a^2(x_0)=0$ then certainly $lim_{x\rightarrow x_0}\a(x)=0$ and taking $\a(x_0)=0$ we get a function continuous at $x_0$.

Let $c(t)$ be a unit speed geodesic with $c(0)=x_0\in N$. We take $\xi$ which in geodesic coordinates around $x_0$ is $\xi=\sqrt{|\lb'-\mu'|}\frac{\partial}{\partial t}$. Note that  $\sqrt{|\lb'-\mu'|}=\pm g(\dot c,\xi)$.   The function  $\phi(t)=\sqrt{|\lb'-\mu'|}\circ c(t) $ if $t\ge 0$ and $\phi(t)=-\sqrt{|\lb'-\mu'|}\circ c(t)$ if $t<0$ is differentiable at $0$ and smooth in a neighborhood of $0$.  Indeed note that  $\phi(t)=g(\dot c(t),\xi(c(t)))$  and  $\phi'(t)=g(\dot c,\n_{\dot c}\xi)=\frac12\a(c(t)$ for $t\ne0$.  Since $\a$ is continuous it follows that $\phi$ is differentiable at $0$ and has continuous derivative.   The function $\a(t)=\a(c(t))$ is differentiable at $0$ and $\a'(0)=0$.   In fact $d\a=-\frac2{n-1}\mu\om$ and  $d\a(\dot c)=-\frac2{n-1}\mu\phi(t)$.  Hence $\a'(t)=-\frac2{n-1}\mu\circ c(t)\phi(t)$.
Note that $\frac12\a(t)=\phi'(t)$. The last equation holds also for $t=0$ since $\a$ is continuous. Hence we get an equation $\phi''(t)=-\frac{\mu(c(t))}{n-1}\phi(t)$ with initial conditions $\phi(0)=0,\phi'(0)=\frac12\a(x_0)$ where the function  $\mu(c(t))$ is real analytic.  If $\a(x_0)=0$ then we get $\phi(t)=0$, a contradiction. Hence $\a(x_0)\ne 0$ and $\xi=\frac1{\a}\n\lb'$ is a smooth vector field in a neighborhood of $x_0$.   Hence $\xi$ is a smooth vector field.$\k$
\medskip
{\it Remark.}  $\Cal{AC}^{\perp}$ manifolds for which $\xi$ is a Killing vector field are described in [J-1]. In the rest of the paper we shall assume that $\dim M=n+1$.
\bigskip
{\bf Theorem 2.} {\it Let us assume that $(M,g)$ is a real analytic, simply connected $\Cal {AC}^{\perp}$ manifold, whose Ricci tensor has two eigenvalues $\lb,\mu$ of multiplicities  $n,1$ respectively and such that $\mu-\frac{2\text{Scal}}{n+3}$ is constant and $\lb-\frac{2\text{Scal}}{n+3}$ is non-constant. Then  $\lb-\mu\ge0$ on the whole of $M$  or $\lb-\mu\le0$ on the whole of $M$.  If $|\lb-\mu|>0$ on $M$ then  $M$ is a warped product  $M=\Bbb R\times_{f^2} M_*$ where $M_*$ is an Einstein manifold.  If $N\ne\emptyset$ then $N=\{x_0\}$ or $N=\{x_0,x_1\}$ and $M=\Bbb R_+\times_{f^2} M_*$ in the first case, while $M=(0,\e)\times_{f^2} M_*$ in the second case, where $M_*=(S^n,g_{can})$ is a sphere with the standard metric of constant sectional curvature. In all  cases the function $f$ satisfies an equation
$$-(n-1)(\frac{f''}f-\frac{(f')^2}{f^2})=\tau f^{-2}-Cf^2\tag 2.4$$ where $C\in \Bbb R-\{0\}$ and $Ric_{M_*}=\tau g_{M_*}$. If $N\ne \emptyset$ and $M_*$ is a sphere with sectional curvature $1$ then the function $f$ satisfies the initial conditions $f(0)=0,f'(0)=1$ in the first case and in the second case $f(0)=0,f'(0)=1, f(\e)=0,f'(\e)=-1$. }

\medskip
{\it Proof.} If $N=\emptyset$ the proof is given in [J-2] (see also [H]). Let us assume first that $\sharp N\ge2$.  Let $x_0\in N$ and $r=d(x_0,N-\{x_0\})=d(x_0,x_1)$ where $x_1\in N$.
We show that the map $exp_{x_0}:B(0,r)\rightarrow M$, where $B(0,r)=\{X\in T_{x_0}M:||X||<r\}$, is a diffeomorphism onto its image.

 If $exp_{x_0}X=exp_{x_0}Y$ then we consider two cases  $X=0,Y\ne 0$ and $X,Y\ne0$.  In the first case the geodesic  $c(t)=exp_{x_0}t\frac Y{||Y||}$ is a loop of length $s=||Y||<r$. In the second case let $c(t)=exp_{x_0}t\frac X{||X||}, c_1(t)=exp_{x_0}t\frac Y{||Y||}$. Let $t_0=||X||,t_1=||Y||$.  Then  $c(t_0)=c_1(t_1)$ and $\dot c(t_0)=\dot c_1(t_1)$ since  $\xi(c(t_0))=\sqrt{|\lb-\mu|}\dot c(t_0)=\sqrt{|\lb-\mu|}\dot c_1(t_0)$. It follows that $c(t_0+s)=c_1(t_1+s)$.  Let us assume that $t_1>t_0$.  Then  $c(0)=x_0=c_1(t_1-t_0)$. Hence $c_1$ is a geodesic loop at $x_0$ of length $t_1-t_0<r$. This contradicts $\xi(c_1(t))=\sqrt{|\lb-\mu|}\dot c_1(t)$ which is valid up to the first value of $t>0$ for which $\lb(c_1(t))-\mu(c_1(t))=0$, if we start from $t=0$  (note that $|\lb(c_1(t))-\mu(c_1(t))|=0$ only if $c_1(t)=x_0$ for $t\in (0,r)$).  Thus  $t_0=t_1$ and $c=c_1$, which means that $X=Y$.  Similarly we obtain $Y=0$ if $X=0$. The mapping  $exp_{x_0}:B(0,r)\rightarrow M$ is a diffeomorphism. In fact, otherwise we would have a conjugate point on a geodesic  $c(t)=exp tX$.  Hence we would have two different geodesics starting at $x_0$ and joining $x_0$ to $x_1=exp sX$, where $s>1, s||X||<r$. It leads to a contradiction as above.

 Note that $\lb-\mu$ depends only on the distance from $x_0$. Hence $exp_{x_0}rS\subset N$ where $S$ is a unit sphere in $T_{x_0}M$. Since $exp_{x_0}$ is a continuous function, $exp_{x_0}rS$ is connected and $N$ consists of isolated points, it follows that $exp_{x_0}rS=\{x_1\}$. Hence $N=\{x_0,x_1\}$ and every unit geodesic $c$ starting from $x_0$ satisfies $c(r)=x_1$. It follows that $exp_{x_0}:B(0,r)\rightarrow M-\{x_1\}$ is a diffeomorphism. In the same way $exp_{x_1}:B(0,r)\rightarrow M-\{x_0\}$ is a diffeomorphism onto $M-\{x_0\}$.

 Note that by the Gauss lemma the leaves of $\1$ are compact spheres $S^n$ via the diffeomorphism $exp_{x_0}:S(0,t)\rightarrow M$ where $S(0,t)=\{X\in T_{x_0}M:\|X\|=t$ where $t\in(0,r)$.  Hence $M=[0,r]\times_{f^2}S^n$.
From [B, Lemma 9.114,  page 269] it follows that on $S^n$ there is a metric of constant sectional curvature $\lb^2g_{can}$ and $f(0)=0=f(r),f'(0)=\frac1{\lb}=-f'(r)$.  Replacing $f$ by a function $\lb f$ we see that we can assume that $f(0)=0,f'(0)=1,f(r)=0,f'(r)=-1$.

Now consider the case $N=\{x_0\}$.  Then  as above one can prove that $exp_{x_0}:T_{x_0}M\rightarrow M$ is a diffeomorphism.  Hence in this case $M=\Bbb R_+\times_{f^2} S^n$ where  $S^n$ has the a canonical metric. The function $f$ satisfies the initial conditions $f(0)=0,f'(0)=1$.$\k$
\medskip
Now we consider the equation  (see [J-2,p.27])
$$-(n-1)(\frac{f''}f-\frac{(f')^2}{f^2})=\tau f^{-2}-Cf^2$$
with initial conditions  $f(0)=0,f'(0)=1$ ; hence we consider the warped product $M=(0,\e)\times_{f}S^{n}$ or  $M=\Bbb R_+\times_{f}S^{n}$  and $\tau=n-1$.

If we write $g=\ln f$ then  $(g)''=-\frac{\tau}{n-1}e^{-2g}+\frac C{n-1}e^{2g}$.  We get
$\frac12((g')^2)'=(\frac{\tau}{2(n-1)}e^{-2g}+\frac C{2(n-1)}e^{2g})'$ and $(g')^2=\frac{\tau}{(n-1)}e^{-2g}+\frac C{(n-1)}e^{2g}+A$
for some constant $A\in\Bbb R$. Consequently, if we assume that $\tau=n-1$  we get an equation for $f$
$$(f')^2=1+Af^2+\frac C{n-1}f^4.$$ Using  homothety we can assume that $C=\pm(n-1)$.   Thus we get
$(f')^2=1+Af^2+f^4$ or $(f')^2=1+Af^2-f^4$.

First let us assume that the quadratic polynomial $1+At+\e t^2$ where $\e\in\{-1,1\}$ has real roots. This happens for $\e=-1$, and for $|A|\ge2$ if $\e=1$.  Hence if $\e=1$ and $A<-2$ we get
$(f')^2=(a^2-f^2)(\frac1{a^2}-f^2)$ for some $a\in(0,1)$.  Let $f$ be the solution of this equation with initial conditions $f(0)=0,f'(0)=1$.

We shall show that $f$ increases until it attains its maximum $a=f(t_0)$ and then decreases, and $f(2t_0)=0,f'(2t_0)=-1$.  In fact $f$ satisfies an equation $f'=\sqrt{(a^2-f^2)(\frac1{a^2}-f^2)}$ if $f<a$. It follows that $\arcsin (\frac fa)'=\sqrt { \frac1{a^2}-f^2}\ge\sqrt{\frac1{a^2}-a^2}$ if $f<a$. Consequently, $\arcsin (\frac fa)\ge t\sqrt{\frac1{a^2}-a^2}$ and $f$ increases until it attains the value $f(t_0)=a$. The solution exists on an interval $[0,t_0]$. Now we define $f(t_0+s)=f(t_0-s)$ for $s\in(0,t_0]$.  Then $f$ is defined on an interval $[0,2t_0]$ and is a smooth function on an open neighborhood of this interval.

Similarly if $\e=-1$ then we get an equation
 $f'=\sqrt{(a^2-f^2)(\frac1{a^2}+f^2)}$ where $a>0$.  Thus $\arcsin (\frac fa)\ge t\sqrt{\frac1{a^2}}$  and again there exists $t_0>0$ such that $f(t_0)=a$. We extend the solution on the interval $[0,2t_0]$  and obtain a solution satisfying  $f(2t_0)=0,f'(2t_0)=-1$.

 We have to show that $f$ is an odd function at $t=0$ and $t=2t_0$.    Note that $f$ satisfies the equation $f''=Af+2\e f^3$.  Since $f$ is real analytic, it suffice to show that $f^{(2k)}(0)=0$ for $k\in\Bbb N$.   We use induction.    We have  $f(0)=0,f''(0)=0$.   Let us assume that  $f^{(2k)}(0)=0$ for $k\le n$ where $n\ge1$.  We will show that $f^{(2n+2)}(0)=0$. Note that  $$f^{(2n+2)}(0)=Af^{(2n)}(0)+2\e\sum_{r=0}^{2n}\sum_{s=0}^rC^{2n}_rC^r_sf^{(s)}f^{(r-s)}f^{(2n-r)},\tag 2.5$$ where $C^n_k=\frac{n!}{(n-k)!k!}$. If $s$ and $r-s$ are odd numbers then $r=r-s+s$ is an even number and $f^{(2n-r)}(0)=0$. Hence in all cases the left hand side of $(2.5)$ equals $0$.  It follows that in both cases considered above we get a compact $\Cal{AC}^{\perp}$ manifold diffeomorphic to $S^{n+1}$.
If $\e=1$ and $A=-2$ then one can easily check that $f(t)=\tanh t$.  In that case we obtain a complete noncompact $\Cal{AC}^{\perp}$ manifold diffeomorphic to $\Bbb R^{n+1}$  (see [J-2]).
If $A>-2$ then we can obtain a solution which tends to infinity which may not be defined on the whole of $\Bbb R_+$.  For example for $A=2$ we get $f(t)=\tan t$.
 Summarizing we have

\medskip
{\bf Theorem 3.}  {\it   On the sphere $S^{n+1}$ there exists a one-parameter family of $\Cal{AC}^{\perp}$ metrics $g_A$ such that $g_A=dt^2+f^2g_{can}$ and $(S^{n+1},g_A)=[0,2t_0]\times_{f^2}S^{n}$ where $f$ satisfies the equation $$(f')^2=1+Af^2+f^4\tag 2.6$$ with $A\in(-\infty,-2)$, or $$(f')^2=1+Af^2-f^4\tag 2.7$$ for $A\in\Bbb R$.
In the case of  equation (2.6) the corresponding $\Cal{AC}^{\perp}$ manifold admits a pair of  Einstein-Weyl structures both conformally Einstein and conformal to the standard sphere with constant sectional curvature.}
\medskip
{Proof.} Note that in the case (2.6) the corresponding $\Cal{AC}^{\perp}$ manifold has a Ricci tensor with eigenvalues $\lb,\mu$ such that $\lb-\mu =Cf^2\ge0$.  Hence from (1.5) it determines a pair of Einstein-Weyl structures  $(g,\om'), (g,-\om')$ where  $\om'=\frac2{\sqrt{n-1}}\om$ and $\om(X)=g(\xi,X)$ was defined in the first part of the paper.  Since $d\om=0$, both Einstein-Weyl structures are conformally Einstein. Let $\om=d\phi$ then the metrics $g_1=\exp(\phi)g, g_2=\exp(-\phi)g$ are Einstein metrics on $S^{n+1}$.  They are conformally equivalent to each other  $g_1=\exp(2\phi)g_2$.  Thus it follows from   [K-R-1], [K-R-2]  that $g_i$ are standard metrics of constant sectional curvature on $S^{n+1}$.

\medskip
{\it Remark.} In the case of equation (2.6) with $A=-2$ we obtain a solution $f(t)=\tanh t$.   The corresponding $\Cal{AC}^{\perp}$ manifold $M=\Bbb R_+\times_{f^2} S^n$ is complete (see [J-2]) and diffeomorphic to $\Bbb R^{n+1}$.  It also admits a pair of Einstein-Weyl structures which are conformally Einstein.     At least one of them has to be non-complete since any complete Einstein manifold which is conformal to a complete Einstein manifold has to be a sphere with standard metric of constant sectional curvature (see [K-R-1]).

\bigskip
\medskip
{\bf Theorem 4.}  {\it   On a manifold $M=\Bbb R\times _{f^2} M_*$ where $(M_*,g_{M_*})$ is an Einstein manifold such that $Ric_{M_*}=\tau g_{M_*}$ with $\tau<0$ there exists a one parameter family of complete $\Cal{AC}^{\perp}$ metrics $g_A$ such that $g_A=dt^2+f^2g_{M_*}$, where $f$ satisfies the equation $$(f')^2=\frac{\tau}{n-1}+Af^2-f^4\tag 2.8$$ with $A\in\Bbb R, A^2>-4\tau,\tau<0$. The eigenvalues $\lb,\mu$  of the Ricci tensor of $(M,g_A)$
satisfy $\lb\le\mu$ hence these manifolds do not admit a pair of Einstein-Weyl structures.}
\medskip
{\it Proof.}  Equation (2.8) is equivalent to the equation   $$f''= Af-2f^3\tag 2.9$$ with appropriate initial conditions. Hence we can apply  [B,Lemma 16.37 p.445 ].   Note that $(f')^2=-(f^2-a^2)(f^2+\frac{\tau'}{a^2})$ if $A^2>-4\tau',A>0,\tau'<0$ where $\tau'=\frac{\tau}{n-1}, a>0$. We can assume that  $a<\frac{\sqrt{-\tau'}}a$. Hence  for $P(f)=-(f^2-a^2)(f^2+\frac{\tau'}{a^2})$ we have $P(a)=P(\frac{\sqrt{-\tau'}}{a})=0$ and $P'(a)P'(\frac{\sqrt{-\tau'}}{a})\ne0$.   Consequently, the equation $(f')^2=P(f)$ admits a periodic solution $f:\Bbb R\rightarrow [a,\frac{\sqrt{-\tau'}}a]$. In that way we get complete manifolds $(M,g_A)$. Note that these manifolds admit compact quotients.$\k$
\medskip
{\it Remark.} Note that an equation of type $(2.9)$ ( for a function $\frac1f$) was studied in [B] and complete and compact manifolds constructed in Theorem 4 are known (as manifolds with $Dr\in C(Q)$ in [B]). However the  $\Cal{AC}^{\perp}$ metrics on spheres $S^{n+1}$ were not found by A. Besse.

\bigskip
\centerline{\bf References.}

\par
\medskip
[B] A. Besse `Einstein manifolds' {\it Springer Verlag} Berlin
Heidelberg (1987)

 \par
\medskip
[G] P. Gauduchon {\it La 1-forme de torsion d'une variete
hermitienne compacte} Math. Ann. {\bf 267} (1984) 495-518.
\par
\medskip
[Go] A. Ghosh {\it Complete Riemannian manifolds admitting a pair of Einstein-Weyl structures} Mathematica Bohemica, 141, 315-325 (2016)
\par
\medskip
[Gr] A. Gray {\it Einstein-like manifolds which are not Einstein}
Geom. Dedicata {\bf 7} (1978) 259-280.
\medskip
[H] S. Hiepko {\it Eine innere Kennzeichnung der verzerrten Produkte }
Math. Ann. {\bf 241}  (1979)  209-215.
\par
\medskip
[J-1] W. Jelonek {\it Killing tensors and Einstein-Weyl geometry}
Colloquium Math., {\bf 81}  (1999) 5-19.
\par
\medskip
[J-2]W.Jelonek  {\it Killing tensors and warped products}
 Ann.  Polon. Math. {\bf 75}   (2000) 15-33.

\par
\medskip
[K-R-1] W. K\"uhnel, H-B. Rademacher {\it Einstein spaces with a conformal group}, R\-e\-sults in Mathematics {\bf 56}  (2009) 421-444.
\par
\medskip
[K-R-2] W. K\"uhnel, H-B. Rademacher {\it Conformal vector fields on pse\-udo-Riema\-nnian
spaces} Differential Geometry and its Applications {\bf 7} (1997) 237-250.

\par
\medskip
[PS1] H. Pedersen and A. Swann {\it Riemannian submersions,
four-manifolds and Einstein-Weyl geometry} Proc.London Math. Soc
(3) {\bf 66} (1993) 381-399.
\par
\medskip
[PS2] H. Pedersen and A. Swann {\it Einstein-Weyl geometry, the
Bach tensor and conformal scalar curvature} J. reine angew.
Math.{\bf 441} (1993) 99-113.

\par
\medskip
[T]  K.P. Tod {\it Compact 3-dimensional Einstein-Weyl structures}
J.  London Math.  Soc.(2) {\bf 45} (1992) 341-351.
\bigskip
 Institute of Mathematics,

 Cracow University of  Technology, Warszawska 24,

31-155 Krak\'ow,POLAND.

E-mail address: wjelon\@pk.edu.pl

\end